\documentclass[11pt]{article}
  \usepackage[ansinew]{inputenc}
  \usepackage[psamsfonts]{amssymb}
 \def\bk{{l\!k}} 

 \def\hC{{\bf C}} 

 \def\bB{I\!\!B} 

 \title{Loop homology algebra of a closed manifold}
 \author{Yves F\'elix, Jean-Claude Thomas and Micheline Vigu\'e-Poirrier}

\begin{document}

 \maketitle

 \begin{abstract} The loop homology of a closed orientable manifold $M$ of
 dimension $d$ is the ordinary homology of the
 free loop space $LM = M^{S^1}$ with degrees shifted by $d$, i.e. $\mathbb
 H_\ast  (LM) = H_{\ast  +d}(LM)$. M. Chas and D. Sullivan
 have defined a loop product on $\mathbb H_\ast  (LM)$ and an intersection
 morphism $I : \mathbb H_\ast  (LM) \to
 H_\ast  (\Omega M)$. The algebra $\mathbb H_\ast  (LM)$ is
commutative   and
 $I$ is a   morphism  of
 algebras. In this paper we produce, for any field  $\bk$, a chain  model
that computes the algebra
 $\mathbb H_\ast  (LM; \bk)$ and the morphism  $I$. We
 show that the kernel of $I$ is a nilpotent ideal and that the image of
$I$ is contained in
 the center of $H_\ast (\Omega M; \bk)$, which is in
 general quite small.
 \end{abstract}

 \vspace{5mm}\noindent {\bf AMS Classification} : 55P35, 54N45,55N33,
17A65,
 81T30, 17B55

 \vspace{2mm}\noindent {\bf Key words} : free loop space, loop homology,
 Hochschild cohomology, Gerstenhaber product.

 \vspace{1cm}
  Let $M$ be a simply connected closed oriented $d$-dimensional
 (smooth) manifold and $\bk$ be a principal ideal domain. In   (\cite{CS}),
Chas
and   Sullivan have constructed a product on
  the desuspension, $\mathbb
 H_\ast  (LM;\bk) = H_{\ast  +d}(LM;\bk)$,  of the
 ordinary homology of the free loop
  space on $M$. This product, called the {\it loop product}, is
   defined at the chain level using both  intersection product on the
chains on $M$ and loop composition. (See 2.2   for more details).
The loop homology of $M$ is the commutative graded algebra $
\mathbb H_\ast(LM; \bk)$. Our first result consists in producing an
explicit model, which allows to compute the loop product (Theorem
4.4) for any field $\bk$.

In their paper Chas and Sullivan   construct also a  morphism of
graded algebras, called the {\it intersection  morphism},
 $$
 I : \mathbb H_\ast (LM;\bk) \to H_\ast (\Omega M;\bk )\,,
$$
 where $H_\ast  (\Omega M;\bk)$ is the usual (non commutative) Pontryagin

algebra. Using a description at the chain level of the
 morphism $I$, we prove

 \vspace{3mm} \noindent {\bf Theorem A.} (Theorem 5.2) {\sl For any field
$\bk$,

a) the kernel of
 $I$ is a nilpotent ideal  of nilpotency index less than or equal to
$d/2$,

b) the image of $I :\mathbb H_\ast  (LM; \bk) \to H_\ast  (\Omega
 M; \bk)$ lies in the center of $ H_\ast  (\Omega M;\bk)$. }

 \vspace{3mm} Theorem A shows that the image of $I$ is in general
 very small comparatively to the expected growth of $H_\ast  (LM;
\bk)$. When $\bk$ is a field of characteristic zero, Theorem A
becomes  more precise. Let us recall that an  element $x \in
\pi_q(M)$ is called a {\it Gottlieb element}
 (\cite{FHT}-p.377), if the map $x \vee id_M : S^q
 \vee M \to M$ extends to the product $S^q \times M$. These elements
generate a   subgroup $G_\ast  (M)$ of $\pi_\ast  (\Omega M)$ via
the isomorphism
 $\pi_\ast  (\Omega M) \cong \pi_{\ast +1}(M)$. Finally, we denote  by
cat $M$ the Lusternik-Schnirelmann category of $M$ normalized so
that cat $S^n$=1.

 \vspace{3mm}\noindent {\bf Theorem B.} (Theorem 6.3) {\it If $\bk$ is a
field of characteristic zero then

a)   the kernel of
 $I$ is a nilpotent ideal  of nilpotency index less than or equal to cat
$M $.

b) (Im $I$) $ \cap\, (\pi_\ast (\Omega M)\otimes {\bk} ) =
G_\ast  ( M)\otimes {\bk}$.

c)  $\displaystyle\sum_{i=0}^n$  dim (Im $  I \cap H_i(\Omega
M;\bk) \leq C n^k  $,   some constant $C>0$ and  $k \leq$ cat $M
$.}

 \vspace{3mm}

Chas and Sullivan have    proved that $I$ is surjective when $M$
is a
 Lie group, and in (\cite{CJY}), Cohen, Jones and Yan have computed
 the algebra $\mathbb H_\ast  (LM;\bk)$ and the  homomorphism $I$
when $M$
 is a sphere or a complex projective space. Using our model we perform
their results, proving in particular:

  \vspace{3mm}\noindent
{\bf  Theorem C.} (Theorem 7.5) {\it  The intersection
  morphism  $I : \mathbb H_\ast  (LM;\mathbb Q) \to H_\ast
(\Omega M;\mathbb Q)$
 is surjective if and only
 if
 $M$ has the rational homotopy type of
    a product of odd dimensional spheres.}

 \vspace{3mm} The starting point   of our work is the Cohen-Jones
 isomorphism
 $$f_*   : \mathbb H_{\ast  }(LM; \bk
)\stackrel{\cong}{\longrightarrow}  H\!H^\ast ({\cal C}^\ast M;
 {\cal C}^\ast M)   $$
 which identifies the loop product
 with the
 Gerstenhaber product
   on $ H\!H^\ast ({\cal C}^\ast M;
 {\cal C}^\ast M)$ (see 1.5)
when   $M$ is  a 1-connected closed oriented manifold of dimension
$d$ (\cite{CJ}).

The Cohen-Jones ring isomorphism    depends a priori on the smooth
structure because its construction  uses in an essential way  the
normal bundle of the  embedding of the manifold into some $\mathbb
R^n$. So even though $H\!H^*({\cal C}^*M, {\cal C}^*M)$ is a
homotopy invariant, there is no evidence that the ring isomorphism
$f_*$ is independent of the smooth structure, and that, for a
homotopy equivalence $g : M \to M'$, the induced map $g_* :
{\mathbb H}_*(LM) \to {\mathbb H}_*(LM')$ is a ring isomorphism.
 Anyway this problem does not affect our results that only depend
on the existence of a ring isomorphism and not on the isomorphism
itself.

Denote by   $\varepsilon :
  {\cal C}^\ast  (M) \to \bk$ the augmentation associated to the inclusion
$\{m_0\}\hookrightarrow M$. Then the following result provides us
 a description of the intersection morphism in terms of
Hochschild cohomology.

  \vspace{3mm}\noindent {\bf Theorem D.} (Theorem 3.5) {\sl There is a
commutative diagram
   of algebras
 $$
 \renewcommand{\arraystretch}{1.6}
 \begin{array}{ccc}
 H\!H^\ast  ({\cal C}^\ast  (M); {\cal C}^\ast (M)) & \stackrel{H\!H^\ast
({\cal C}^\ast (M);
 \varepsilon)}{\longrightarrow} &
 H\!H^{\ast   }({\cal C}^\ast  (M);\bk)\\
 \cong \uparrow {\scriptstyle f_*} && \uparrow \cong \\
 \mathbb H_\ast  (LM; \bk) & \stackrel{I}{\to} & H_{\ast   }(\Omega M;
\bk)\,,
 \end{array}\,.
 \renewcommand{\arraystretch}{1}
 $$
  }

  To prove Theorem A we also use the following algebraic
result concerning the center of the enveloping algebra of a graded
Lie algebra.

\vspace{3mm}\noindent {\bf Theorem E.} (Theorem 6.2)  {\sl Let $L$
be a finite type graded Lie algebra defined on a field of
characteristic zero, then the center of $UL$ is contained in the
enveloping algebra on the radical of $L$.}

\vspace{2mm} The paper is  organized as follows:

\noindent  1 - Hochschild  cohomology and Gerstenhaber product.

\noindent 2 - Loop product and Hochschild cohomology.

\noindent 3 - The intersection   morphism  $I :  \mathbb H_\ast
(LM) \to  H_\ast  (\Omega M)$.

\noindent 4 - A chain model for computing the loop product and
$I$.

\noindent 5 - The kernel and the image of $I$.

\noindent 6 -  Determination of $I$  when $\bk$ is a field of
characteristic zero.

 \noindent 7 - Examples and applications.

 \noindent 8 - Hochschild cohomology and Poincar\'e duality.

\vspace{3mm} We would like to thank the referee for helpful
comments that have considerably improve   the paper.

 \vspace{5mm}

{\bf 1 - Hochschild  cohomology and Gerstenhaber product }

 In this section we  fix some  notations  and recall the standard
 definitions of Hochschild  cohomology and of Gerstenhaber product.

\vspace{3mm}
 \noindent{\bf 1.1.}   Let $\bk $   be a principal ideal domain; modules,
tensor product, linear homomorphism,...  are defined over $\bk$.
For notational simplicity, we avoid to mention $\bk$. If $V $ is a
lower or upper graded module ($V_{i}=V^{-i}$) the
 suspension $s$ is defined by  $(sV)_n =V_{n+1}\,, \,\, (sV)^n
 =V^{n-1}$.

 \vspace{3mm}
\noindent{\bf 1.2.} Let  $(A, d)$ be a differential graded
augmented
 cochain algebra and $(N,
 d) $ be a
 differential graded $ A $-bimodule, $
 A =   \{A^i\}_{i\geq 0},   N =\{N^j\}_{j \in {\mathbb Z}}$
and $\bar A =$ ker$(\varepsilon : A \to \bk )$. The {\it two-sided
normalized bar
 construction},
$$
 \overline{\bB} (N; A; N) =  N \otimes T (s  \bar A ) \otimes N\,,
\quad \overline{\bB} _k (N; A; N) =  N \otimes T^k (s  \bar A )
\otimes N\,,
 $$
is defined as follows: For $k\geq 1$, a generic element    $
 m[a_1|a_2|...|a_k]n \in
 \overline{\bB} _k (N; A; N)
 $  has  degree    $| m|+ | n|+ \sum_{i=
 1}^k (|s a_i|)\, $.  If $k=0$, we write $m[\,]n = m\otimes 1\otimes n \in
N \otimes T^0s\bar A\otimes N$. The
 differential $d= d_0+d_1$ is defined by:
 $$
 d_0  : \overline{\bB}_k  (N;A;N) \to \overline{\bB}_k (N;A;N)\,,
 \quad d_1  : \overline{\bB}_k (N;A;N) \to \overline{\bB}_{k-1}
 (N;A;N)\,,
 $$
 with
 $$
 \renewcommand{\arraystretch}{1.6}
 \begin{array}{rll}
 d_0  ( m[a_1|a_2|...|a_k]n)& = d( m) [a_1|a_2|...|a_k]n  -
\displaystyle\sum
 _{i=1}^k (-1)^{\epsilon _i}  m[a_1|a_2|...|d(
 a_i)|...|a_k]n\\ &+ (-1) ^{\epsilon _{k+1}}  m[a_1|a_2|...|a_k]d(n)
 \\[3mm]
 d_1  (m[a_1|a_2|...|a_k]n) &= (-1) ^{| m|}  ma_1[a_2|...|a_k]n +
 \displaystyle\sum _{i=2}^k (-1) ^{\epsilon _i}
 m[a_1|a_2|...|a_{i-1}a_i|... | a_k]n \\
 &- (-1)^{\epsilon _{k}}  m[a_1|a_2|...|a_{k-1}]
 a_k n
 \end{array}
 \renewcommand{\arraystretch}{1}
 $$
  Here
 $\epsilon _i = | m| + \sum _{j<i} (|s a_j|)$.

  \vspace{3mm}
\noindent{\bf 1.3.} For any differential graded algebra $ A $, let
$A^{op}$ be the opposite
 graded
 algebra, $a\cdot^{op}b = (-1)^{\vert a\vert . \vert b\vert} b\cdot a$,
and
 $A^e = A \otimes A^{op}$ be the enveloping algebra. Any differential
graded $A$-bimodule $N$ is a differential graded $A^e$-module.

 Let $A$ and $N$ as in 1.2. The {\it Hochschild  cochain complex}
$\hC^\ast  (A;N)$ of $A$ with coefficients
 in  $N$ is  the differential graded module   (\cite{Ge},\cite{Vi}):
 $$
  {\hC}^\ast  (A;N) =  \mbox{Hom} _{A^e}  (\overline{\bB} (A; A;
 A), N) \,, \quad \hC^n(A,M) = \prod_{p-q=n}
 \mbox{Hom}_{A^e}(\overline{\bB}(A,A,A)^p,N^{q})\,,$$
equipped with the standard differential $D$ defined by $Df = d
\circ f - (-1) ^{|f|} f\circ d $. The cohomology of the complex
$\hC^\ast  (A;N)$ is called the {\it Hochschild cohomology} of $A$
with values in $N$
, and is denoted ${H\!H}^\ast (A;N) $.

This definition extends the
classical one since:

\vspace{3mm} \noindent{\bf 1.4. Lemma} (\cite{FHT-1}-Lemma 4.3)
{\it  If $A$ is a differential graded algebra such that $A$ is a
$\bk$-free graded module  then the multiplication in $A$ extends
in  a semi-free resolution  of $A^e$-modules
$$
m :\overline{\bB}(A,A,A) \longrightarrow A\,.
$$}

This means that $m$ is a quasi-isomorphism of differential graded
$A$-bimodules which  well behaves with quasi-isomorphisms of
differential graded $A$-bimodules. In particular, we have the
following lifting lemma

\vspace{3mm}\noindent {\bf Lifting Homotopy Lemma.} {\sl For any
quasi-isomorphism $\varphi: A' \to A$ there exists a unique (up to
homotopy in the category of differential graded bimodules)
quasi-isomorphism $ \hat m :\overline{\bB}(A,A,A) \longrightarrow
A'$ such that $ m \simeq \varphi \circ \hat m $.}

 \vspace{3mm}

\noindent{\bf 1.5.}  Recall  that $\overline{\bB}(A) =
 \overline{\bB}(\bk;A;\bk)  := \left( T(s\bar A), d\right) $  is a
differential graded coalgebra  with
 $$ d([a_1|a_2|\cdots |a_k]) = - \sum
 _{i=1}^k (-1)^{\epsilon _i}  [a_1|a_2|...|d( a_i)|...|a_k] + \sum
 _{i=2}^k (-1) ^{\epsilon _i}  [a_1|a_2|...|a_{i-1}a_i|... |
 a_k]\,.$$
The canonical isomorphism of graded modules $ \mbox{Hom} _{A^e}
(\overline{\bB} (A; A; A), N)  = \mbox{Hom} (T(s\bar A), N)\,, $
carries on $\mbox{Hom} (T(s\bar A), N)$ a differential $D'$.
Observe that the differential $D'$ is not the canonical
differential $D$  of $\mbox{Hom} (\overline{\bB} (A), N)$ except
when $N$ is the trivial bimodule $\bk$.

 If $N=A$, Gerstenhaber \cite{Ge} has proved that the  usual cup
product on  $\mbox{Hom}   (T(s\bar A), A)$ makes $(\mbox{Hom}
(T(s\bar A), A),D')$   a differential graded algebra such that
   the induced product on ${H\!H}^\ast (A;A) $,
called the {\it Gerstenhaber product}, is commutative.

\vspace{5mm}

{\bf 2 - Loop product and Hochschild cohomology}

In this section, for the reader convenience and   in order to
precise notations, we recall the definition of the loop product
defined by Chas and Sullivan, the interpretation given by Cohen
and Jones  and the relation between the loop product and the
Gerstenhaber product.

 Let $X$ be topological space  and ${\cal C}_\ast
(X)$ be
 the  singular chain coalgebra (coefficients in  $\bk$) with coproduct
$\Delta _X$:
$$
{\cal C}_k  (X) \stackrel{\mbox{\small diagonal}}{\longrightarrow}
{\cal C}_k  (X\times X)  \stackrel{AW}{\longrightarrow}
 \renewcommand{\arraystretch}{0.6}
\begin{array}[t]{c}
\oplus\\
{\scriptstyle k_1+k_2=k}
\end{array}
\renewcommand{\arraystretch}{1}
{\cal C}_{k_1} (X)\otimes {\cal C}_{k_2}  (X) \,, \quad c \mapsto
\sum_i c_i \otimes c'_i\,,
$$
where $ AW $ denotes the usual Alexander-Whitney chain
equivalence.

The singular cochain algebra on $X$ is the
differential graded algebra  ${\cal C}^\ast  (X):= \left( {\cal
C}_\ast  (X)\right)^\vee$ with product ({\it cup product}),
$$
\cup : {\cal C}^p   (X) \otimes {\cal C}^q   (X) {\longrightarrow}
{\cal C}^{p+q}  (X) \,, \quad (f\cup g)  (c) = \sum _i (-1)^{\vert
c_i\vert \cdot \vert g\vert} f(c_i) g(c'_i)\,.
$$
Observe that the choice of a point $x_0$ in $X$ determines an
augmentation $\varepsilon : {\cal C}^\ast (X) \to \bk$.

\vspace{3mm} \noindent{\bf 2.1.}  Let $M$ be a   1-connected
$\bk$-Poincar\'e duality space of  formal  dimension $d$, i.e.,
$M$ is
1-connected with a class $[M] \in H _d (M) $ such that cap product with
$[M]$ induces  an isomorphism
$$
-\cap [M] : H^\ast (M) \to H_{d-\ast} (M) \,.
$$
(We restrict  to 1-connected spaces in order to simplify the definition.
In this case the formal dimension $d$ does not depend on
$\bk$.) Every simply-connected closed oriented $d$-dimensional manifold is
a $\bk$-Poincar\'e duality space.

Let us denote also by $[M]$ a cycle in ${\cal C}_d( X)$ representing the
class $[M] \in H_d (M)$, and
 denote by ${\cal P} : H _\ast
(M)
\to H^{d-\ast }(M)$ the linear isomorphism, inverse of
$-\cap [M] $.

 The
intersection product of two homology classes  of $M$ is nowadays
defined from the cup product in the following way:
$$
\renewcommand{\arraystretch}{1.6}
\begin{array}{ccc}
H_k(M) \otimes H_l(M) &\stackrel{- \bullet -}{\longrightarrow}&
H_{k+l-d}(M) \\
\hspace{-8mm}{\scriptstyle {\cal P}\otimes {\cal P } }\downarrow
&&
\uparrow {\scriptstyle -\cap [M]}\\
H^{d-k} (M) \otimes H^{d-l} (M) &\stackrel{- \cup
-}{\longrightarrow}& H^{2d-(k+l)}(M)
\end{array}
\renewcommand{\arraystretch}{1}
$$

The intersection product makes ${\mathbb H}_\ast (M) = H_{\ast
+d}(M)$ a commutative graded algebra with unity $[M] \in {\mathbb
H}_0 (M)$. Originally, the intersection product was defined for a
polyhedron.  In their paper (\cite{CJ}) Cohen and Jones uses the
following point of view. Let $i : M \hookrightarrow {\mathbb
R}^{n+d}$ be the embedding of $M$ into a codimension $n$ euclidean
space. We denote by $\nu^n \to M$ the normal bundle of this
embedding, and by $Th(\nu^n )$ the associated Thom space, which is
Spanier-Whitehead dual to $M_+$, i.e. $M$ with a disjoint base
point. More precisely, denote by $M^{-TM}$ the spectrum
$$M^{-TM} = \Sigma^{-(n+d)} Th (\nu^n )\,.$$
Then the diagonal map $\Delta : M \to M\times M$ induces a map of
spectra $\Delta^* : M^{-TM} \land M^{-TM} \to M^{-TM}$ that makes
$M^{-TM}$ into a ring spectrum. When $M$ is orientable, the ring
structure is compatible, via the Thom  isomorphism, with the
intersection product:
$$
\begin{array}{ccc}
H_q(M^{-TM}) \otimes  H_r(M^{-TM}) & \to & H_{q+r} (M^{-TM})\\
\downarrow {\scriptstyle \tau \otimes \tau} & & \downarrow
{\scriptstyle
\tau}\\
H_{q+d} (M) \otimes H_{r+d}(M) & \stackrel{\mbox{\scriptsize
intersection}}{\longrightarrow} & H_{q+r+d}(M)\,.
\end{array}
$$

\vspace{3mm} \noindent{\bf 2.2.} The  {\it loop product}  has been
defined geometrically by Chas and Sullivan  (\cite{CS} in the
following way. Consider again the evaluation map  $\mbox{\small
ev}  :
 LM \to M$   and let $\sigma :
 \triangle ^n \to LM$ and $\tau : \triangle ^m \to LM$ be
singular simplices. Then $\mbox{\small ev}\circ\sigma$ and
$\mbox{\small ev} \circ \tau$ are singular simplices of $M$. At
each point $(s,t) \in \Delta^n \times \Delta^m$ where $q\circ
 \sigma (s) = q\circ \tau (t)$, the composition of the loops $\sigma (s)$
and
 $\tau (t)$ can be defined. If we assume
that the map  $(\mbox{\small ev}\circ\sigma, \mbox{\small ev}
\circ \tau) : \triangle ^n \times \triangle ^m \to M\times M $ is
transverse to the diagonal map  $M \to M\times M$ then
 (as shown by Chas and Sullivan) this defines a
chain $\sigma \cdot \tau\in
 C_{n+m-d}(LM)$,  and produces a
 commutative and associative multiplication
 $$
\mathbb H_q(LM) \otimes \mathbb H_r(LM) \to
 \mathbb H_{q+r}(LM)\,, \hspace{1cm} a\otimes b\mapsto a\bullet
 b\,,
 $$
  whose unity is the image of $[M]$ via the
homomorphism induced by the canonical section.

A homological presentation of the loop product in the smooth case
goes as follows. Let $N \to M$ be the normal bundle to the
inclusion $\Delta : M \to M \times M$ and $(ev)^*(N)$ the pull
back bundle on $LM\times_MLM$. We write $N_D, N_S, (ev)^*(N)_D,
(ev)^*(N)_S$ for the corresponding disk and sphere bundles. We
have homotopy equivalences
$$\varphi_1 : (N_D, N_S) \to (V, \partial V)$$
$$\varphi_2 : ((ev)^*(N)_D, (ev)^*(N)_S) \to ((ev)^{-1}(V),
(ev)^{-1}(\partial V))\,,$$ where $V$ is a tubular neighborhood of
$\Delta (M)$ into $M\times M$. The map $\varphi_1$ is the usual
exponential map and $\varphi_2$ is defined by
$$\varphi_2(x,v,c) = \gamma (x,v)^{-1} \circ c \circ \gamma (x,v)\,,$$
with $x \in \Delta (M), v\in N_x, c\in LM\times_MLM$ and where
$\gamma (x,v)$ denotes the geodesic ray of length $\parallel
v\parallel$ starting from $x$ with tangent vector $v$. Then the
Chas-Sullivan loop product can be defined as the composite
$$
\renewcommand{\arraystretch}{1.6}
\begin{array}{l}
H_q(LM) \otimes H_r(LM) \to H_{q+r}(LM\times LM, (ev)^{-1}
(M\times M
{\scriptstyle\backslash} \Delta (M) ))\\
\stackrel{\cong (1)}{\longrightarrow} H_{q+r}((ev)^{-1}(V),
(ev)^{-1}(\partial V)) \cong H_{q+r} ((ev)^*(N)_D,
(ev)^*(N)_S) \\
\stackrel{\cong (2)}{\longrightarrow} H_{q+r-d}(LM\times_MLM)
\stackrel{\mbox{\scriptsize composition}}{\longrightarrow}
H_{q+r-d} (LM)\,,
\end{array}
\renewcommand{\arraystretch}{1}
$$
where $(1)$ is usual excision and $(2)$ the Thom isomorphism for
an oriented fiber bundle.

\vspace{3mm} \noindent{\bf 2.3.} In (\cite{CJ}), Cohen and Jones
give another description of the Chas-Sullivan product. Since we
need   relative versions of their construction we recall here the
main steps of their construction. Using the notations of 2.1, let
$Th(ev^*(\nu^n))$ be the Thom space of the pull back bundle
$ev^*(\nu^n) \to LM$ and define the Thom spectrum
$$LM^{-TM} = \Sigma^{-(n+d)}Th (ev^*(\nu^n))\,.$$

\vspace{2mm}\noindent {\bf 2.4. Theorem.} (\cite{CJ}, Theorem 1.3)
{\sl $LM^{-TM}$ is a ring spectrum whose multiplication
$$\mu : LM^{-TM} \land LM^{-TM} \to LM^{-TM}$$
is compatible with the Chas-Sullivan product in the sense that the
following diagram commutes
$$
\begin{array}{ccc}
H_q(LM^{-TM}) \otimes H_r(LM^{-TM}) &
\stackrel{\mu^*}{\longrightarrow} &
H_{q+r} (LM^{-TM})\\
\downarrow {\scriptstyle \tau} & & \downarrow {\scriptstyle \tau}\\
H_{q+d}(LM) \otimes H_{r+d}(LM) &
\stackrel{\bullet}{\longrightarrow} & H_{q+r+d}(LM)\,,
\end{array}
$$
where $\tau$ denotes the Thom isomorphism.}

\vspace{3mm}\noindent {\bf 2.5. Theorem.} (\cite{CJ}, Theorem 3)
{\sl $LM^{-TM}$ is the geometric realization of a cosimplicial
spectrum $({\mathbb L}_M)_*$. The $k$-simplices are given by
$$({\mathbb L}_M)_k = (M^k)_+\land M^{-TM}$$
and the ring structure is realized on the cosimplicial level by
pairings
$$\mu_k : ((M^k)_+\land M^{-TM}) \land ((M^r)_+\land M^{-TM}) \to
(M^{k+r})_+ \land M^{-TM}$$ defined by
$$\mu_k ((x_1, \ldots , x_k;a) \land (y_1, \ldots , y_r; b)) = (x_1,
\ldots , x_k,y_1, \ldots , y_r; \Delta^*(a\land b))$$ where
$\Delta^*$ is the ring structure defined on $M^{-TM}$   in 2.1.}

\vspace{5mm}\noindent {\bf 2.6.}  Let $\Delta^k$ be the standard
$k$-simplex
$$\Delta^k = \{ (x_1, \ldots , x_k)\, \vert \, 0 \leq x_1 \leq x_2
\leq\ldots \leq x_k \leq 1\, \}$$ and consider the maps
$$f_k : \Delta^k \times LM \to M^{k+1}\,, \hspace{5mm}
 f_k((x_1, \ldots , x_k),c) = (c(0), c(x_1), \ldots , c(x_k))\,.$$
Denote by $\tilde f_k : LM \to \mbox{Map} (\Delta^k, M^{k+1})$ the
adjoint of $f_k$. The product of the $\tilde f_k$ induces an
homeomorphism
$$f : LM \to \mbox{Map}_{\Delta_*}(\Delta^*, M^{*+1})
\hspace{4mm}(\cite{CJ})\,,$$ which induces, when $M$ is simply
connected, a linear isomorphism in homology
$${\mathbb H}_*(LM) \stackrel{\cong}{\longrightarrow} H\!H^*({\cal
C}^*(M), {\cal C}^*(M)) \hspace{3mm} \cite{CJ}\,.$$ Pulling back
the normal bundle $\nu^n \to M$ along the maps $ev$ and $p_1$ in
the diagram
$$
\begin{array}{ccc}
\Delta^k \times LM & \stackrel{f_k}{\longrightarrow} & M^{k+1}\\
\downarrow {\scriptstyle ev} && \mbox{} \hspace{2mm}\downarrow
{\scriptstyle p_1}\\
M & = & M\,,
\end{array}$$
we get maps of Thom  spectra
$$f_k : (\Delta^k)_+ \land LM^{-TM} \to M^{-TM} \land (M^k)_+\,.$$
By gluing together the adjoint maps, we get a map of spectra
$$f : LM^{-TM} \to \prod_k \mbox{Map} ((\Delta^k)_+ , M^{-TM}\land
(M^k)_+)\,.$$ Now, using the homotopy equivalence  ${\cal
C}_*(M^{-TM}) \simeq {\cal C}^{-*} (M_+)$, and taking singular
chains we get a morphism connecting  ${\cal C}_*(LM^{-TM})$ to the
Hochschild cochains on ${\cal C}^*(M)$,
$$f_* : {\cal C}_{*-k} (LM^{-TM}) \to   {\hC}^k({\cal C}^*(M), {\cal
C}^*(M))\,.$$ The commutativity of  $f_*$   with the differentials
is not a trivial point. One reason is that, on the spectrum level,
the Atiyah duality map $M^{-TM} \to \mbox{Map}(M, S^0)$ is a map
of ring spectra and of bimodules of spectra (\cite{Co}).

\vspace{3mm}\noindent {\bf 2.7. Theorem.}  (\cite{CJ}, Theorem 12)
{\sl Let $M$ be a compact simply connected manifold. Then the
morphism $f_*$ is a quasi-isomorphism and the following diagram of
chain complexes is commutative
$$
\begin{array}{ccc}
{\cal C}_*(LM^{-TM}) \otimes {\cal C}_*(LM^{-TM}) &
\stackrel{\mu_*}{\longrightarrow} & {\cal C}_*(LM^{-TM})\\
{\scriptstyle f_*\otimes f_*} \downarrow {\scriptstyle \cong} & &
{\scriptstyle f_*} \downarrow {\scriptstyle \cong}\\
{\hC}^*({\cal C}^*(M), {\cal C}^*(M))\otimes {\hC}^*({\cal
C}^*(M), {\cal C}^*(M))& \stackrel{\cup}{\longrightarrow}&
{\hC}^*({\cal C}^*(M), {\cal C}^*(M))
\end{array}
$$}

\vspace{5mm} {\bf 3 - The intersection   morphism  $I :  \mathbb
H_\ast  (LM) \to
 H_\ast  (\Omega M)$}

 \vspace{3mm}
\noindent {\bf 3.1.} Let $i : N \hookrightarrow M$ be the
injection of an open set in $M$. We define the mapping space
$L_NM$ as the pullback
$$
\begin{array}{ccc}
L_NM & \stackrel{i'}{\longrightarrow} & LM\\
\downarrow & & \downarrow {\scriptstyle ev} \\
N & \stackrel{i}{\longrightarrow} & M
\end{array}
$$
The space $L_NM$ is the space of loops that originate in $N$. By
restriction, the  loop product induces a product on $ {\mathbb H}_*(L_NM)$ so
that the morphism $H_*(i') : {\mathbb H}_*(L_NM) \to {\mathbb H}_*(LM)$
becomes a
multiplicative morphism.

We will now follow verbatim the lines of the Cohen-Jones
construction in order to compute ${\mathbb H}_*(L_NM)$ in terms of
Hochschild cohomology. With the notation of section 2.3, we define
the Thom ring spectra
$$L_NM^{-TM} = \Sigma^{-(n+d)} Th(i'^*ev^*(\nu^n))\,.$$ The Thom map
$\tau$ is then a multiplicative isomorphism
$$ H_*(\tau )  : H_*(L_NM^{-TM}) \stackrel{\cong}{\to} H_{*+d}(L_NM)\,.$$
The morphisms $f_k$ restrict naturally to morphisms
$$f_k : L_NM^{-TM} \land \Delta^k \to N^{-TM} \land (M^k)_+$$ that induce,
exactly in the same way as in the original case $N=M$, a
multiplicative isomorphism
$$H_*(L_NM^{-TM}) \stackrel{\cong}{\to} H\!H^*({\cal C}^*(M), {\cal
C}_*(N^{-TM}))\,.$$ Now by the Spanier-Whitehead duality, there is
a homotopy equivalence
$${\cal C}_*(N^{-TM}) \simeq {\cal C}^*(M, M{\scriptstyle \backslash}
N))\,,$$ and the map $f_*$ restricts naturally to a
quasi-isomorphism $$ f_*^N : {\cal C}_{*-k} (L_NM^{-TM}) \to
{\hC}^k({\cal C}^*(M), {\cal C}^*(M, M{\scriptstyle\backslash} N
))\,.$$
 We therefore have

\vspace{3mm}\noindent {\bf 3.2. Theorem.} {\sl When $M$ and $N$
are simply connected, there is a commutative diagram of algebras
in which horizontal lines are isomorphisms
$$
\renewcommand{\arraystretch}{1.6}
\begin{array}{cccccc}
\varphi_N : & {\mathbb H}_*(L_NM) &
\stackrel{\cong}{\longleftarrow} & H_*(L_NM^{-TM}) &
\stackrel{\cong}{\longrightarrow} & H\!H^*({\cal
C}^*(M), {\cal C}^*(M, M{\scriptstyle \backslash} N))\\
& \downarrow {\scriptstyle H_*(i'}) && \downarrow && \downarrow\\
\varphi_M: & {\mathbb H}_*(LM) & \stackrel{\cong}{\longleftarrow}
& H_*(LM^{-TM}) & \stackrel{\cong}{\longrightarrow} & H\!H^*({\cal
C}^*(M), {\cal C}^*(M)) \,,
\end{array}
\renewcommand{\arraystretch}{1}
$$
where, as usual, ${\mathbb H}_p = H_{p+d}$.}

\vspace{3mm}\noindent In the same way, working in the relative
case, we get

\vspace{3mm}\noindent {\bf 3.3. Theorem.} {\sl There is a sequence
of isomorphisms of algebras
$$
\renewcommand{\arraystretch}{1.6}
\begin{array}{ll}
\varphi_{M,N} : {\mathbb H}_*(LM,L_NM) \stackrel{\tau
}{\longleftarrow}& H_*(LM^{-TM}, L_NM^{-TM}) \\&
\stackrel{\cong}{\to} H\!H^*( {\cal C}^*(M), {\cal
C}^*(M{\scriptstyle \backslash} N))
\end{array}
\renewcommand{\arraystretch}{1}
$$
making commutative the diagram
$$
\begin{array}{ccc}
  {\mathbb H}_*(L_NM) & \stackrel{\varphi_N}{\longrightarrow} &
H\!H^*({\cal C}^*(M),
{\cal C}^*(M, M{\scriptstyle \backslash} N))\\
  \downarrow  &&   \downarrow\\
  {\mathbb H}_*(LM) & \stackrel{\varphi_M}{\longrightarrow} & H\!H^*({\cal
C}^*(M), {\cal
C}^*(M)) \\
\downarrow && \downarrow \\
  {\mathbb H}_*(LM, L_NM) & \stackrel{\varphi_{M,N}}{\longrightarrow} &
H\!H^*({\cal C}^*(M),
{\cal C}^*( M{\scriptstyle \backslash} N))\\
\end{array}
$$}

\vspace{3mm}\noindent{\bf 3.4.} In \cite{CS},  Chas and Sullivan
define the
 intersection morphism $I :  \mathbb H_\ast  (LM) \to H_\ast
(\Omega M)$ by
 associating to an $q$-cycle in $LM$ its intersection with the space
of based loops at $m_0$, $\Omega (M,m_0)=\Omega M$, viewed as a
codimension $d$-submanifold. It follows directly from the
definition that
   $ I $
 transforms the loop product into the  Pontryagin product.

 A slightly different  exposition of the intersection morphism
 works as follows.  Fix a Riemannian metric on $M$, choose   a
 geodesic disc $D^d$ centered at the base point $m_0$ and consider the
 map
 $$\mu :  D^d \times \Omega M \to LM\,, \hspace{1cm} \mu (x,\omega ) =
 \gamma_x^{-1} \ast
 (\omega \ast \gamma_x)\,,$$ where $\gamma_x$ denotes the geodesic ray
 from $m_0$ to $x$, and $\ast  $ denotes the composition of paths. Let us
denote
 by $E$ the subspace  $L_{M\backslash \{ m_0\}}M$. We consider the
commutative
 diagram of complexes
 $$
 \renewcommand{\arraystretch}{1.6}
 \begin{array}{ccccc}
 {\cal C}_\ast (D^d) \otimes {\cal C}_\ast  (\Omega M) &
 \stackrel{EZ}{\longrightarrow} & {\cal C}_\ast  (D^d\times \Omega M)
 & \stackrel{{\cal C}_\ast  (\mu)}{\longrightarrow} & {\cal C}_\ast
(LM)\\
\alpha \otimes id \downarrow &&&& \downarrow {   \alpha}\\
 {\cal C}_\ast (D^d,\partial D^d) \otimes {\cal C}_\ast  (\Omega M) &&
 \stackrel{\psi}{\longrightarrow} && {\cal C}_\ast  (LM, E)
 \end{array}
 \renewcommand{\arraystretch}{1}
 $$ where $EZ$ means the Eilenberg-Zilber map, $\alpha$ the canonical
surjections, and  $\psi$ is the quotient map. Clearly  $H_\ast
(\psi) : H_\ast (\Omega M) \to \mathbb
 H_\ast  (LM,E)$ is an isomorphism  of algebras. The
 intersection morphism   $I$ coincides with the composition
$H(\psi)^{-1}\circ H(\alpha
)$:
 $$
  H_q(LM) \stackrel{H(\alpha
 )}{\longrightarrow} H_q(LM, E)
 \stackrel{H(\psi)}{\longleftarrow} \left(H_\ast (D^d, \partial
D^d)\otimes H_\ast(\Omega M)\right) _q\cong
  H_{q-d}(\Omega M)\,.
$$

 Our next result  describes the homomorphism   $I$ in terms of
 Hochschild cohomology. Denote by   $\varepsilon :
 {\cal C}^\ast (M) \to \bk$ the augmentation induced by the inclusion
$\{m_0\} \hookrightarrow M$. We then have:

 \vspace{3mm}
\noindent {\bf 3.5 Theorem.} {\sl Let $M$ be a simply connected
closed oriented $d$-dimensional
 manifold.   There exists an isomorphism of graded algebras $\Theta$  that
makes commutative the
following diagram
 $$
 \renewcommand{\arraystretch}{1.6}
 \begin{array}{ccc}
 H\!H^\ast  ({\cal C}^\ast  (M); {\cal C}^\ast (M)) &
  \stackrel{H\!H^\ast  ({\cal C}^\ast  (M); \varepsilon)}{\longrightarrow}
& H\!H^{\ast }({\cal C}^\ast  (M);\bk)\\
 \cong \uparrow {\scriptstyle f_*} && {\scriptstyle \Theta}\uparrow
\cong \\
 \mathbb H_\ast  (LM) & \stackrel{I}{\to} & H_{\ast   }(\Omega M)\,.
 \end{array}
 \renewcommand{\arraystretch}{1}
 $$
  }

 \vspace{3mm}\noindent {\bf Proof.}  We take $N = M \backslash \{ m_0\}$
in Theorem 3.3. We remark that the map $H\!H^\ast  ({\cal C}^\ast
(M); \varepsilon)$ is the composite
$$H\!H^*({\cal C}^*(M), {\cal
C}^*(M)) \to H\!H^*({\cal C}^*(M), {\cal C}^*( M{\scriptstyle
\backslash} N)) \stackrel{\bar\varepsilon}{\to} H\!H^*({\cal
C}^*(M), \bk)\,,$$ where $\bar\varepsilon$ is an isomorphism. We
define the isomorphism $\Theta$ to be the composition
$$
H_*(\Omega M) \stackrel{H_*(\psi)}{\longrightarrow} {\mathbb
H}_*(LM,E) \stackrel{\varphi_{M, M\backslash \{m_0\} }
}{\longrightarrow} H\!H^*({\cal C}^*(M),{\cal C}^*(M{\scriptstyle
\backslash} N)) \stackrel{\bar\varepsilon}{\longrightarrow}
H\!H^*({\cal C}^*(M),\bk)\,.
$$
With these definitions, the above diagram commutes trivially.
 \hfill $\square$

 \vspace{5mm}

{\bf 4 - A chain model for computing the loop product and $I$.}

In this section we construct, for any field of coefficients $\bk$,
an explicit model for the loop product at the chain level.

\vspace{3mm} \noindent{\bf 4.1.}   Recall the  Adams Cobar
construction
 $\Omega C$ on a coaugmented differential graded coalgebra $C = \bk \oplus
\bar C$. This  is
 the differential
 graded algebra $(T(s^{-1}\bar C),d)$, where $d=d_1+d_2$ is  the unique
 derivation determined by:
 $$d_1s^{-1}c = -s^{-1}dc\,, \mbox{
 and }
 d_2s^{-1}c = \sum_i (-1)^{|c_i|}s^{-1}c_i \otimes s^{-1}c_i'\,,
\hspace{3mm}
 c\in \bar C\,,$$
 where the reduced coproduct  of $c\in \bar C$ is written $\bar
 \Delta c = \sum_i c_i\otimes c_i'$.  For sake of simplicity we
 put  $\langle x_1|x_2|\cdots |x_n\rangle := s^{-1}x_1 \otimes \cdots
\otimes
 s^{-1}x_n\,.$

 \vspace{3mm}
\noindent{\bf 4.2.} Assume $\bk $ is a field, and $M$ is a
$1$-connected compact $d$-dimensional manifold.  Denote   by $f :
(T(V),d)
 \to {\cal C}^\ast  (M)$   a free minimal model for the singular cochain
 algebra on $M$ (\cite{FHT-1}), i.e. $(T(V),d)$ is a differential
 graded algebra, $f$ is a quasi-isomorphism of differential
 graded algebras, and $d(V) \subset T^{\geq 2}(V)$. The
 differential graded algebra $(T(V),d)$ is uniquely defined, up to
isomorphism,  by the
 above properties. Moreover, $V^p \cong  H^{p-1}(\Omega M)$,
 (\cite{FHT-1}). Denote by  $S
$   a complement of the vector space generated by the cocycles of
degree $d$. The differential graded ideal $J = (T(V))^{>d} \oplus
S$  is acyclic and the quotient algebra $A
 = T(V)/J$ is a finite dimensional graded differential algebra.

  \vspace{3mm}
\noindent{\bf 4.3.} Since $A$ is finite dimensional, the dual
algebra $A^\vee $ is a differential graded coalgebra  and we
consider the
 differential graded algebra
$$
\Omega A^{\vee} =(T(W),d)\,,
$$
with in  particular,
  $W \cong \mbox{Hom}(s\bar A,\bk)$, and $\Omega A^{\vee}
=\mbox{Hom}(\overline{\bB}A,\bk) = (\mbox{Hom}( T(s\overline A),
\bk), D)$
  (see 1.5).

   We choose a homogeneous linear  basis $e_i$ for
 $\overline{A}$, and its  dual basis $w_i$ for $W$. This determines
the constants of structure $\alpha_{ij}^k$ and $\rho_i^j$:
$$
\renewcommand{\arraystretch}{1.3}
\begin{array}{rl}
\langle w_i, se_k\rangle &= - (-1)^{\vert
 w_i\vert}\delta_{ik}\,, \quad e_i \cdot e_j = \sum_k \alpha_{ij}^k e_k\,,
\quad  d(e_i)
 = \sum_j \rho_i^je_j\\
 d(w_i)& = \sum_{jk}a_i^{jk}w_jw_k + \sum_j\beta_i^j w_j\,, \quad
 a_i^{jk} = (-1)^{\vert e_j\vert +\vert e_je_k\vert}
 \alpha_{jk}^i\,, \quad \beta_i^j = (-1)^{\vert w_j\vert}\rho_j^i\,.
\end{array}
\renewcommand{\arraystretch}{1}
$$

 \vspace{3mm}\noindent {\bf 4.4 Theorem.} {\it Let $\bk$ be a field and
$M$ be a 1-connected closed oriented manifold of dimension $d$.
With notation introduced above

 a) the derivation  $D$  uniquely defined on the tensor product of graded
algebras $A \otimes T(W)$   by
 $$\left\{
 \renewcommand{\arraystretch}{1.6}
 \begin{array}{l}
 D(a\otimes 1) = d(a) \otimes 1 + \sum_j (-1)^{\vert a\vert +
 \vert e_j\vert} [a,e_j]\otimes w_j\,, \hspace{1cm} a\in A\,,\\
  D(1\otimes b) = 1 \otimes d(b) - \sum_j (-1)^{\vert e_j\vert}e_j
 \otimes [w_j,b]\,, \hspace{1cm} b \in TW\,,
 \end{array}
 \right.
 \renewcommand{\arraystretch}{1}
 $$
is a differential. Here  $[\,\,,\,\,]$ denotes the Lie bracket in
the graded algebras  $A$ and  $T(W)$.

b) the graded algebra  $H_\ast (A\otimes T(W), D)$ is isomorphic
to the loop algebra   $\mathbb H_\ast (LM)$.}

 \vspace{2mm}\noindent {\bf Proof.}

a) is proved by a direct but laborious computation.

b) is a direct consequence of theorem 2.7. \hfill {$\square$}

Observe that this model is dual to those constructed by one of us,
\cite{Vi}.

 \vspace{3mm}\noindent {\bf 4.5 Proposition.}  {\it Let $\bk$ be a field
and  $M$ be  a 1-connected closed oriented manifold of dimension
$d$.
 There is a cohomology spectral
 sequence of graded algebras such that
 $$E_2 = H\!H^\ast  (H^\ast  (M),H^\ast  (M)) \Rightarrow \mathbb
 H_\ast  (LM)\,.$$}
\noindent {\bf Proof.} The spectral sequence is obtained by
 filtering   the complex $(\mbox{Hom}\, (T(s\overline A,A),D')$ by the
differential ideals $\mbox{Hom}\,
 (T^{\leq  p}(s \overline A),A)$  (see 1.5). Since $H^\ast  (A) = H^\ast
(M)$, it follows that $E_1 = \mbox{Hom}
 (\overline{\bB}(H^\ast  (M)),\bk)\otimes
 H^\ast  (M) $ and $ E_2 = H^\ast  (\mbox{Hom}\,
 (\overline{\bB}(H^\ast  (M)),\bk)\otimes
 H^\ast  (M) ,D)\,.
$ \hfill $\square$

 \vspace{3mm}
\noindent {\bf 4.6  Example.} If  $M$ is a formal space, (for
instance $ M $ is   a simply connected compact
 K\"ahler manifold for $\bk = \mathbb
 Q$, \cite{DGMS}), one can choose
 $A = H^\ast  (M)$ and thus  the algebras
 $H\!H^\ast  (H^\ast  (M); H^\ast  (M))$ and $\mathbb H_* (LM)$
are isomorphic graded vector spaces. If we put $H^*=H^*(M)$ and
$H_*=H_*(M)$ the loop algebra $\mathbb H_*(LM)$ is isomorphic to
the graded algebra $H(H^* \otimes T(s\overline H_*),D)$ with $
D(a\otimes 1) =0$, $a\in H^*$ and $D(1\otimes b)= -\sum_J
(-i)^{|e_j|} e_j \otimes [w_j, b]$, $ b \in \overline H_*$.

 \vspace{3mm}
\noindent {\bf 4.7.   The commutative case.} Suppose  that the
algebra ${\cal C}^\ast (M)$ is
 connected by a sequence of quasi-isomor\-phisms to a commutative
 differential graded algebra $(A,d)$. This is the case if either $\bk$ is
of
 characteristic zero, or else if  $\bk$ is a field of characteristic
 $p > d$ (\cite{A}, Proposition 8.7). We can also suppose that $A$ is
finite dimensional, $A^0 = \bk$, $A^1 =
 0$, $A^{>d} = 0$ and $A^d = \bk \omega$. Then
 formulas of Theorem 4.4 simplify as:
 $$
 \renewcommand{\arraystretch}{1.6}
 \left\{
 \begin{array}{l}
 D(a\otimes 1) = d(a) \otimes 1\,,\\
 D(1\otimes b) = 1 \otimes d(b) - \sum_j (-1)^{e_j} e_j \otimes
 [w_j,b]\,.
 \end{array}
 \right.
 \renewcommand{\arraystretch}{1}
 $$

\vspace{3mm} We can now interpret the intersection morphism in
terms of models.

\vspace{3mm}   \vspace{3mm}\noindent  {\bf 4.8 Theorem.} {\sl Let
$\bk$ be a field and $M$ be a 1-connected closed oriented manifold
of dimension $d$. There is a commutative diagram of algebras
 $$\begin{array}{ccc}
 \renewcommand{\arraystretch}{1.6}
 \mathbb H_\ast  (LM) & \stackrel{\cong}{\longrightarrow} &
 H_\ast  (A\otimes T(W),D)\\
 {\scriptstyle I}\downarrow & & \downarrow {\scriptstyle
 H(\varepsilon_A \otimes 1)}\\
 H_\ast  (\Omega M) &\stackrel{\cong}{\longrightarrow} & H_\ast
(T(W),d)\,.
 \end{array}
 \renewcommand{\arraystretch}{1}
 $$
 }

 \vspace{2mm}\noindent {\bf Proof.}  Recall that
 Hochschild cohomology $H\!H^\ast  (A; M)$ is covariant in $M$ and
 contravariant in $A$. Moreover, if $f : A \to B$ is a
 quasi-isomorphism of differential graded algebras and $g : M \to
 M'$ is a quasi-isomorphism of $A$-bimodules, we have isomorphisms
 $$H\!H^\ast  (B;M) \stackrel{\cong}{\longrightarrow} H\!H^\ast  (A;M)
 \stackrel{\cong}{\longrightarrow} H\!H^\ast  (A;M')\,.
$$
Starting with
 Theorem 3.5, we obtain the following commutative diagram
 $$
 \begin{array}{ccccccc}
 \mathbb H_\ast  (LM) & \stackrel{\cong}{\longrightarrow} &
 H\!H^\ast  ({\cal C}^\ast  (M);{\cal C}^\ast  (M)) &
 \stackrel{\cong}{\longrightarrow} & H\!H^\ast (A;A) &
 \stackrel{\cong}{\longrightarrow} & H_\ast  (A  \otimes T(W),D)\\
 {\scriptstyle I}\downarrow && \downarrow {\scriptstyle
 H\!H^\ast  ({\cal C}^\ast  (M),\varepsilon  )} &   & \downarrow
{\scriptstyle
 H\!H^\ast  (A, \varepsilon_A )} & & \downarrow {\scriptstyle
H(\varepsilon_A
 \otimes 1)}\\
 H_\ast  (\Omega M) & \stackrel{\cong}{\longrightarrow} & H\!H^\ast
({\cal
 C}^\ast  (M);\bk) & \stackrel{\cong}{\longrightarrow} & H\!H^\ast
(A;\bk) &
 \stackrel{\cong}{\longrightarrow} &   H_\ast  (T(W),d)
\end{array}
\,.$$ \hfill $\square$

\vspace{5mm}

{\bf  5 - The kernel and the image of $I$}

 \vspace{3mm}
\noindent{\bf 5.1.}  If $J$ is an
 ideal of an algebra $A$, we put $J^1=J $ and $J^{n+1}=J J^{n}\,, n\geq 1$
and, in the case
 $J$ is nilpotent, we define  $\mbox{Nil}\, (J)  =
 \mbox{sup}\,\{ n\,\vert \, J^n \neq 0\,\}$.

  \vspace{3mm}
\noindent {\bf 5.2 Theorem.} {\it Let $\bk$ be a field and $M$ be
a simply connected closed   oriented $d$-dimensional
 manifold.

a) The kernel of
 the intersection  morphism $I$ is nilpotent and $\mbox{Nil}\,
(\mbox{Ker}\,
 I) \leq d/2$.

b)  The image of $I$ is
 contained in the center of $H_\ast  (\Omega M)$.}

  \vspace{3mm}\noindent
 {\bf Proof.} a)  By Theorem 4.8, the kernel of $I$ is generated by the
classes of
  cocycles
 in $\bar A  \otimes T(W)$. Since $A^1 = 0$ and $A^{>d} = 0$, the
 nilpotency of the kernel of $I$
   is less than or equal to $d/2$.

 b)  Let $e_i$ and $w_i$ be the elements defined in 4.3 and  $[\alpha ]$
be an element in the image of $H(\epsilon _A\otimes id)$. Then
$\alpha$ is a cocycle in   $T(W)$ and there exist elements
 $\alpha_i$ in $T(W)$ such that
 $\bar\alpha = 1 \otimes \alpha + \sum_i e_i \otimes \alpha_i$ is a cycle
 in $A
 \otimes T(W)$. A short calculation shows that the
 component of $e_i$ in $d(\bar\alpha )$ is
 $$(-1)^{\vert e_i\vert} \left( d(\alpha_i) - [w_i,\alpha ] +
 \sum_j\beta_i^j\alpha_j + \sum_{j,k}a_i^{j,k}(-1)^{\vert u\vert \vert
 w_k\vert}
 \alpha_jw_k + \sum_{j,k}a_i^{kj}(-1)^{\vert w_k\vert}
 w_k\alpha_j\,\right)\,.$$
Since this component must be $0$, by Lemma 5.3 below there exists a
surjective  morphism
 $$
 H(T(W),d) \otimes \bk[u] \to H(T(W),d)
$$
 that maps $u$ to $[\alpha]$. This implies that $[\alpha ]$ is in the
 center of $H(T(W),d) \cong H_\ast  (\Omega M)$. \hfill $\square$

 \vspace{3mm}
\noindent {\bf 5.3 Lemma } {\it Assume $\bk$ is a field.  Let
$\alpha$ be a cycle in
 $(T(W),d)$ and let $u$ be a variable in the same
 degree. Then with the notations of 4.3  and 4.4, we have:
 \begin{enumerate}
 \item There exists a surjective quasi-isomorphism
 $$\varphi : (T(w_i, u, w_i'),D) \to (T(W),d) \otimes (\bk[u],0)\,,
 \hspace{1cm}  \vert w_i'\vert = \vert u\vert
 +
 \vert w_i\vert +1 \,,$$
 such that    $\varphi (u) = u$, $\varphi (w_i) = w_i$ and $\varphi (w_i')
=
 0$, and with $D$ defined by
 $$D(w_i') = [w_i,u] - \sum_j\beta_i^jw_j' -
 \sum_{j,k}a_i^{j,k}(-1)^{\vert u\vert \vert w_k\vert} w_j'w_k -
 \sum_{j,k}a_i^{kj}(-1)^{\vert w_k\vert} w_kw_j'\,.$$
 \item There exists a  morphism  of differential graded algebras
 $\rho : (T(w_i ,u , w_i'),D) \to (T(W),d)$  such that $\rho(u)= \alpha$
and
 $\rho(w_i) = w_i$  if and only if there are elements $\alpha_i
 \in T(W)$ satisfying
 $$d(\alpha_i) = [w_i,\alpha ] - \sum_j\beta_i^j\alpha_j -
 \sum_{j,k}a_i^{j,k}(-1)^{\vert u\vert \vert w_k\vert} \alpha_jw_k
 - \sum_{j,k}a_i^{kj}(-1)^{\vert w_k\vert} w_k\alpha_j\,.$$
 \end{enumerate}}

 \vspace{2mm}
\noindent {\bf Proof.} We define $D(w_i')$ by the above
 formula. Proving that $D^2 = 0$ is an
 easy and standard computation. The  morphism
 $$\varphi : (T(w_i, u,w_i'),D) \to (T(w_i),d)\otimes (\bk[u],0)$$
 defined by $\varphi (w_i) = w_i$, $\varphi (u) = u$ and $\varphi
 (w_i') = 0$ is a surjective  homomorphism  of differential graded
 algebras. To prove that $\varphi$ is a quasi-isomorphism, we
 filter each differential graded algebra  by putting $u$ in
 filtration degree $0$ and the other variables in filtration degree one.
We
 are then reduced  to prove that
 $$
\bar\varphi : (T(w_i, u, w_i'),D) \to (T(w_i),0) \otimes
(\bk[u],0)\,,
 \hspace{1cm} d(w_i) = 0, d(w_i') = [w_i, u]\,,
$$
 is a quasi-isomorphism. Denote by $K$ the kernel of
 $\bar\varphi$ and consider the   short exact sequence of complexes
 $$
0 \to (K \otimes E,D) \to (T(w_i, u, w_i')\otimes
 E,D)\stackrel{\bar\varphi \otimes 1}{\longrightarrow}
 ((T(w_i) \otimes \bk[u] )\otimes E,\bar D)\to 0\,,
$$ where $E$ is the
 linear span of the elements $1, sw_i, su$ and $sw_i'$,
 and where $D$ is defined by
 $$D(sw_i) = w_i \otimes 1\,, D(su) = u\otimes 1\,, D(sw_i') = w_i'
 - (-1)^{\vert w_i\vert}w_i\otimes su + (-1)^{\vert u\vert \vert
 w_i'\vert + \vert u\vert} u\otimes sw_i\,.$$ By construction,
 $(T(w_i, u, w_i')\otimes E,D)$ and    $(T(w_i)
 \otimes \bk[u] \otimes E,\bar D)$ are contractible and therefore
 quasi-isomorphic. Now a non-zero cocycle of lowest degree in $K$
 remains a non-trivial cocycle in the complex $(K \otimes E,D)$.
 Therefore $H_\ast  (K) = 0$ and $\varphi$ is a quasi-isomorphism.

 Part (2) of the Lemma follows directly from the expression of $D$.
 \hfill $\square$

\vspace{5mm}

  {\bf 6 -  Determination of $I$  when $\bk$ is a field of characteristic
zero}

In this section $\bk$ is a field of characteristic zero.

\vspace{3mm} \noindent{\bf 6.1.} By Theorem 5.2,  the   image of
$I$ is contained in the center of
 $H_\ast  (\Omega M)$. On the other hand, by the Milnor-Moore theorem (e.g
\cite{FHT}-Theorem 21.5),
 $H_\ast  (\Omega M) $ is the universal enveloping algebra of the
homotopy Lie
 algebra $L_M = \pi_\ast  (\Omega M)\otimes \bk$ (\cite{FHT}-p.294).

Let $L$ be any graded algebra. The
 center, $Z(UL)$, of the  universal enveloping algebra $UL$ contains the
 universal enveloping algebra of the center of the Lie algebra, $UZ(L)$.
 However the inclusion can be strict. Consider for instance the Lie
 algebra $L = {\mathbb L}(a,b)/([b,b],[a,[a,b]])$, with $\vert a
 \vert = 2$ and $\vert b\vert = 1$. The element $(ab-ba)b$ is in
 the center of $UL$, but not in $UZ(L)$. We denote by $R(L)$ the sum of
all solvable ideals in $L$, (\cite{FHT}-p.495).

 \vspace{3mm}
\noindent {\bf 6.2 Theorem.} {\sl If  $L=\{L_i\}_{i \geq 1}$ is  a
graded
  Lie algebra over  a field of characteristic zero satisfying
$\dim L_i <\infty $  then $Z(UL) \subset UR(L)$.}

 \vspace{2mm}\noindent {\bf Proof.} It is well known   that in
 characteristic zero, $UL$ decomposes into a direct sum
$$UL =
 \renewcommand{\arraystretch}{0.5}
\begin{array}[t]{c}
\oplus\\
{\scriptstyle {k\geq 0}}
\end{array}
\renewcommand{\arraystretch}{1}
\,\, \Gamma^k(L)
$$
 where the $\Gamma^k(L)$ are
 sub-vector spaces that are stable for the adjoint representation
 of $L$ on $UL$:  $\Gamma^0 (L) = \bk$, $\Gamma^1(L) = L$, and
 $\Gamma^n(L)$ is the sub-vector space generated by the elements
 $\varphi (x_1, \ldots, x_n) = \sum_{\sigma \in \Sigma_n }
 \varepsilon_\sigma x_{\sigma (1)}\cdots x_{\sigma (n)}, \,\,\,\,
 x_i \in L$. The coproduct $\Delta$ of $UL$ respects the
 decomposition, i.e.
 $$\Delta : \Gamma^n(L) \to
 \renewcommand{\arraystretch}{0.5}
\begin{array}[t]{c}
\oplus\\
{\scriptstyle p+q = n}
\end{array}
\renewcommand{\arraystretch}{1}
 \Gamma^p(L) \otimes
 \Gamma^q(L)\,.$$
 If we denote by $\Delta_p$ the component of $\Delta$ in $  \Gamma^p(L)
\otimes \Gamma^{n-p}(L)$ then
 $$
\Delta_p (\varphi (x_1, \ldots, x_n)) = \sum_{\tau \in Sh_p}
 \varepsilon_{\tau}\,
 \varphi (x_{\tau (1)},\ldots , x_{\tau (p)}) \otimes
 \varphi (x_{\tau (p+1)},\ldots , x_{\tau (n)})\,,
$$
where $Sh_p $ denotes
 the set of $p$-shuffles
 of the set $\{ 1, 2, \ldots , n\}$.
This implies that the
 composition
 $\Gamma^n(L) \stackrel{\Delta_p}{\longrightarrow} \Gamma^p(L)
 \otimes \Gamma^{n-p}(L) \stackrel{\rm\scriptstyle
 multiplication}{\longrightarrow} UL$ is the multiplication by
 $
 \left(
 \renewcommand{\arraystretch}{0.7}
 \begin{array}{l}n\\p
 \end{array}
 \renewcommand{\arraystretch}{1}
 \right)
 $.
We then consider the composite $c$
 $$
 c : \Gamma^n(L) \stackrel{\Delta_1}{\longrightarrow} L \otimes
 \Gamma^{n-1}(L) \stackrel{1 \otimes \Delta_1}{\longrightarrow} L
 \otimes L \otimes \Gamma^{n-2}(L) \to \cdots \to L^{\otimes
 ^n}\,.
$$

 Let  $\alpha\in UL$  be an element in the center of $UL$,
 $\alpha = \sum_{i=1}^n \alpha_i$ with $\alpha_i \in \Gamma^i(L)$. Since
   $\Gamma^i(L)$ is stable by adjunction, each $\alpha_i$ is in
 the center of $UL$. Therefore we can assume that $\alpha \in
 \Gamma^n(L)$. We write $c(\alpha)$ as a sum of monomials $ x_{i_1}\otimes
\ldots \otimes x_{i_n}$.
 Since ${\rm mult}\circ c: \Gamma^n(L) \to UL$ is the multiplication by
 $n!$, the element $\alpha$ belongs to the Lie algebra generated by
 the $x_{i_j}$.
 Suppose that in the decomposition of $c(\alpha )$ the
 number of monomials  is minimal, then for each $r$, $1\leq
 r\leq n$, the elements
 $x_{i_1} \otimes \ldots \otimes x_{i_{r-1}} \otimes
 x_{i_{r+1}}\ldots \otimes x_{i_n}$
 are linearly independent. Since $[\alpha, x] = 0$,  $x \in L$, we obtain
the equation:
 $$ 0 = \sum_{k=1}^n\left(\,\,\sum_i \,
 (-1)^{\vert x\vert \cdot (\vert x_{i_1}\vert + \ldots +\vert
 x_{i_{k-1}}\vert )}
 x_{i_1}\otimes \ldots \otimes
 [x,x_{i_k}]\otimes \ldots \otimes x_{i_n}\,\,\right)\,.$$
 Let us assume that the $x_{i_k}$are ordered by increasing  degrees then
the elements $x_{i_k}$ with maximal degree belong  to $Z(L)$. The
above
 equation
 shows also that $[x_{i_k},x] $ belongs to the subvector space generated
 by the elements $x_{i_l}$ with higher degree. A decreasing induction on
the degree  shows
 that all the $x_{i_k}$ belong to $R(L)$.  \hfill $\square$

 \vspace{3mm}
\noindent {\bf 6.3.} Denote by $X_0$ the $0$-localization of a
simply connected space $X$. The Lusternik-Schnirelmann category of
$X_0$, cat\,$X_0$,  is less than or equal to  the
Lusternik-Schnirelmann of $X$, cat\, $X$. Moreover the invariant
cat\,$X_0$ is easier to compute than cat\,$X$,
(\cite{FHT}-$\S$-27).

\vspace{3mm} \noindent  {\bf  Theorem.}  {\it  Let $M$ be a simply
 connected oriented  closed manifold and
$\bk$ is a field of characteristic zero then

a)   The kernel of $I$ is a nilpotent ideal and Nil( Ker ($I$))
$\leq$ cat $M_0$ .

b) (Im $I$) $ \cap\, (\pi_\ast (\Omega M)\otimes {\bk} ) =  G_\ast
( M)\otimes {\bk}$.

c)  $\displaystyle\sum_{i=0}^n$  dim (\, Im$  I \cap H_i(\Omega
M;\bk) \,) \leq C n^k  $, some constant $C>0$ and  $k \leq$ cat
$M_0$.}

\vspace{3mm}\noindent {\bf Proof.} a) By (\cite{FHT}-Theorems 29.1
and 28.5), ${\cal C}^\ast  (M;\mathbb Q)$ is connected
  by a sequence of quasi-isomorphisms to a connected finite dimensional
  commutative differential graded algebra $(A,d)$ satisfying
  Nil $(\overline A) \leq n $ for $n> \mbox{cat}\, M_0$. Thus we conclude
as in the proof of theorem 5.2).

 b) The differential graded algebra $\Omega
 (A^{\vee}) = (T(W),d)$ is the universal enveloping algebra on the graded
Lie
 algebra ${\cal L}_M = ({\mathbb L}(W),d)$, and the differential
 graded algebra $(T(W \oplus \bk u\oplus  W'),D)$ is the universal
enveloping
 algebra of the differential graded Lie algebra ${\cal L}^1_M =
 (\mathbb L (W \oplus ku\oplus W'),D)$, (e.g \cite{FHT}-p.289), with
 $$ \left\{
 \renewcommand{\arraystretch}{1.6}
 \begin{array}{l}
 d(w_i) = \sum_j \beta_i^j w_j + \sum_{j,k} \frac{1}{2} a_i^{jk}
[w_j,w_k]\,,\\
   D(w_i') = [w_i,u] - \sum_j \beta_i^j w_j' -\sum_{j,k} a_i^{kj}
 (-1)^{\vert w_k\vert } [w_k, w_j']\,.
   \end{array}
  \renewcommand{\arraystretch}{1}
   \right.
   $$
By construction ${\cal L}_M$ is a free Lie  model for $M$
 and   ${\cal L}^1_M$ is a free Lie  model
 for $M \times S^n$ with $n = \vert u\vert +1$, (\cite{FHT}-$\S$24).
Moreover there exists a bijection between homotopy classes of
maps:
 $$[X\times S^n,X] \cong [(\mathbb L (W \oplus \bk
 u\oplus W'),D), (\mathbb L (W),d)] \,.$$
 Therefore a  homomorphism  $\varphi :
 ({\mathbb L}(W \oplus \bk u\oplus  W'),D) \to ({\mathbb L} (W),d)$
 such that  $\varphi (u) = \alpha$ and $\varphi (w)=w$, $w\in W$,
corresponds
 to a  map $f : M \times S^n \to M$ which extends $id_M\vee g :
 M\times S^n \to M$, such that $[g] =\alpha$
 modulo the identifications   $\pi_n(M)\otimes \bk \cong \pi_{n-1}(\Omega
 M)\otimes \bk \cong H_{n-1}({\mathbb L}(W),d)$. This means
 exactly that
 ${\rm Image}\, I \cap (\pi_\ast  (\Omega M)\otimes \bk) =
  { G}_\ast  ( M) \otimes \bk$.

 c) By, Theorems 36.4, 36.5 and 35.10 of \cite{FHT} we know that if
$L=\pi_\ast  (\Omega M)\otimes \bk $
 then $R(L)$ is finite dimensional and $\mbox{ dim}\,
 R(L)_{\rm \scriptstyle even} \leq \mbox{\rm cat} M_0$. We conclude using
the graded Poincar\'e-Birkhoff-Witt theorem (\cite{FHT}-Theorem
21.1): $Z(UL) \subset UR(L) \cong \Lambda (R(L)_{\mbox{\scriptsize
odd}})\otimes \bk[(R(L)_{\mbox{\scriptsize even}}]$.
 \hfill $\square$

\vspace{5mm}

  {\bf 7 - Examples and applications}

In this section we assume that $\bk$ is a field.

 \vspace{3mm}
\noindent {\bf 7.1  The spheres $S^n$.}  The loop homology of
spheres has been computed by Cohen,
 Jones and Yan using a spectral sequence similar to those  obtained  in
4.5. Our
 model leads  to
   a direct computation of the loop homology and the intersection
homomorphism  $I$.

   Since the differential graded algebra ${\cal C}^\ast  (S^n)$ is
 quasi-isomorphic to $(H^\ast  (S^n), 0) =
   (\land u/u^2,0)$, $\vert u\vert = n$, by Example 4.6, $\mathbb
   H_\ast  (LS^n)$ is isomorphic as an algebra to
 $$H^\ast  (\land  u\otimes T(v), D)\,,  \,\, \vert v\vert = n-1\,, \,\,
\vert u\vert =-n\,, \,\, D(u) = 0\,,\,\, D(v) =
   u\otimes [v,v]\,.$$
 When $n$ is odd, $D = 0$,
   $\mathbb H_\ast  (LS^n) \cong \land u\otimes
   T(v)$ and $ I = \varepsilon \otimes 1 : \land u\otimes
   T(v) \to T(v)$. When $n$ is even, $D(v^{2n}) = 0$, $D(v^{2n+1}) =
2u\otimes
   v^{2n+2}$.
   Therefore a set of generators of $\mathbb H_\ast  (LS^n)$ is
   given by the elements
    $c = 1 \otimes v^2\,, b = u\otimes v\,, a = u\otimes 1$ and,
    $$\mathbb H_\ast  (LS^n) \cong \land (b)\otimes \bk [a,c]/
    (2ac, a^2, ab)\,, \,\,  |a| = -n\,, \vert b\vert =
 -1\,, \vert c\vert = 2n-2\,.$$
    The  homomorphism  $I : \mathbb H_\ast  (LS^n) \to H_\ast
(\Omega S^n)
    = T(v)$ is given by: $I(c) = v^2$, $I(a)=I(b)= 0$.

   \vspace{3mm}
\noindent {\bf 7.2 An example  where $I$ is the trivial
homomorphism.} Let $M$ be the connected sum  $M =(S^3\times
S^3\times
   S^3)\#(S^3\times S^3\times S^3)$. The wedge $N = (S^3\times S^3\times
   S^3)\vee  (S^3\times S^3\times S^3)$  is then obtained by attaching a
 $9$-dimensional cell to $M$ along the homotopy
 class determined by the collar between the two components of $M$. Recall
that
 $$\pi_\ast  (\Omega N)\otimes {\mathbb Q} \cong     Ab (a,b,c)
   \coprod
   Ab (e,f,g)\,,$$ where $Ab(u,v,w)$ means the abelian Lie algebra
generated
 by $u$, $v$ and $w$ considered in degree 2.
 The inclusion $i : M \to N$ induces a surjective map
 $\pi_\ast  (\Omega M)\otimes \mathbb Q \to \pi_\ast  (\Omega N)\otimes
 \mathbb Q$,  This means that
   the attachment of the cell is inert in the sense of \cite{FHT}-p.503.
Therefore,
 (\cite{FHT}-Theorem 38.5),
   $$\pi_\ast  (\Omega M)\otimes \mathbb Q \cong Ab (a,b,c) \coprod
   Ab (e,f,g) \coprod {\mathbb L}(x)$$
   with $\vert x\vert = 7$.
   In particular  $R(L)$ is zero, and by Theorems 5.2 and 6.2,
when $\bk$ is of characteristic zero,  the  homomorphism $I$ is
trivial.

 \vspace{3mm}
\noindent {\bf 7.3 Product of two manifolds.} If $A$ and $B$ are
differential graded algebras, we have a natural
 isomorphism of algebras
 $$H\!H^\ast  (A\otimes B, A \otimes B) \cong H\!H^\ast  (A,A) \otimes
 H\!H^\ast  (B,B)\,.$$ Therefore, from Theorem 2.7, for closed oriented
simply
 connected manifolds
 $M$ and $N$, the isomorphism
 $$\mathbb H_\ast  (L(M\times N) ) \cong \mathbb H_\ast  (LM) \otimes
 \mathbb H_\ast  (LN )$$
 respects the Chas-Sullivan product. Moreover the intersection
homomorphism
 $I_{M\times N}$ identifies to $I_M\otimes I_N$.

\vspace{3mm}\noindent {\bf 7.4. Lie groups.}  Let $\bk$ be a field
of characteristic zero and $G$ be a connected  Lie group. Since
$G$ has the rational homotopy type of a
   product of odd dimensional spheres, we obtain
   $$\mathbb H_\ast  (LG;\mathbb Q) \cong \land (u_1, \ldots , u_n)
   \otimes \mathbb T(v_1, \ldots , v_n)\,,$$
 and $I_G$ is onto. This example generalizes in:

   \vspace{3mm}\noindent
{\bf 7.5 Theorem.} {\it  Let $\bk$ be a field of characteristic
zero and $M$ be a simply connected closed   oriented
$d$-dimensional
 manifold. The intersection
    homomorphism  $I : \mathbb H_\ast  (LM) \to H_\ast (\Omega M)$
 is surjective if and only
 if
 $M$ has the rational homotopy type of
    a product of odd dimensional spheres.}

   \vspace{2mm} \noindent {\bf Proof.} When $M$ has the rational homotopy
 type of the product of odd dimensional spheres,
   then $I$ is clearly surjective. Conversely, if $I$ is
   surjective, then $\pi_\ast  (\Omega M)\otimes {\mathbb Q} = {G}_\ast
(   M)\otimes \mathbb Q$. Thus, $\pi_\ast  (M)\otimes \mathbb Q =
G_{\rm
 odd}\otimes \mathbb Q\,,$ (\cite{FHT}, Proposition 29.8). Let $\{ f_i :
 S^{n_i}
 \to M$,
 $ i=1, \cdots , r\}$
 represent  a given linear   basis of $\pi_\ast  (M)\otimes \mathbb Q$,
and let
 $\varphi_i : S^{n_i} \times M \to M$
   be maps that restrict  to $f_{i}\vee id_M$ on
   $S^{n_i}\vee M$. Then
   the composition
   $$S^{n_1} \times \ldots \times S^{n_r} \hookrightarrow S^{n_1} \times
\ldots
   \times
   S^{n_r}\times M \stackrel{1 \times \varphi_r}{\longrightarrow} S^{n_1}
 \times \ldots \times
   S^{n_{r-1}}\times M \stackrel{1 \times
   \varphi_{r-1}}{\longrightarrow}\ldots \stackrel{1 \times
   \varphi_1}{\longrightarrow}M$$
 induces an isomorphism on the homotopy groups. Therefore, $M$ has the
 rational homotopy type of a product of odd dimensional spheres.
 \hfill $\square$

\vspace{5mm}
 {\bf 8 - Hochschild cohomology and Poincar\'e duality}

When two $A$-bimodules $M$ and $N$ are quasi-isomorphic as
bimodules, then the Hochschild cohomologies $H\!H^*(A;M)$ and
$H\!H^*(A;N)$ are isomorphic.  In this section we  relate the
Hochschild cohomology of the singular cochains algebra on $X$ with
coefficients in itself and   with coefficients in the singular
chains on $X$  when $X$ is a Poincar\'e duality space. The usual
cap product with the fundamental class is not a bimodule morphism.
However the cohomology $  H\!H^n  ({\cal C}^\ast  (M); {\cal
C}_\ast (M))$ and $H\!H^{n  -d}({\cal C}^\ast  (M); {\cal
 C}^\ast  (M)) $ are shown to be isomorphic. This point is not directly
related to the constructions built in sections 4 to 8, but has
its own importance.  For this reason  we have
added this point in the last section of this paper.

\vspace{3mm}

 \noindent{\bf 8.1.} Let $V$ be a   graded module, then
   $V^{\vee}$ denotes the graded  dual,
 $V^{\vee} = \mbox{Hom}_{\bk}(V,\bk)$, and $\langle - ;-  \rangle :
 V^\vee\otimes V \to \bk$ denotes the duality pairing. We denote by
$\lambda _V : V \to V^{\vee \vee}$ the natural inclusion defined
by $\langle \lambda_V ( v), \xi\rangle = (-1) ^{|\xi|} \langle
\xi, v \rangle $.

\vspace{3mm}
 \noindent{\bf 8.2.} Let $X$ be topological space.
The  ${\cal C}^\ast  (X)$-bimodule  structures on  ${\cal C}_\ast
(X)$ and ${\cal
 C}^\ast  (X)^{\vee} $ are explicitly defined by:
 $$
\renewcommand{\arraystretch}{1.4}
\begin{array}{l}
 f\cdot c \cdot g  := (-1)^{|c|(|f|+|g|)+|f|+ |f|\,|g|}(g\otimes
id\otimes f) ( \Delta_X \otimes id ) \circ \Delta _X
(c)\,,\hspace{5mm}   c\in {\cal
C}_\ast  (X) \,,\\
\langle f\cdot \alpha \cdot g ; h\rangle   := (-1)^{|f|}\langle
\alpha ; g\cup
 h\cup f\rangle\,, \hspace{5mm}  f,g,h \in {\cal C}^\ast  (X), \alpha
 \in {\cal C}^\ast  (X)^{\vee}\,.
\end{array}
\renewcommand{\arraystretch}{1}
$$
 Remark that the associativity properties of  $AW$ and of
$\Delta_X$ imply directly that ${\cal C}_\ast  (X)$ is a  graded
${\cal C}^\ast  (X)$-bimodule.

Let $1\in {\cal C}^0(X)$ be the $0$-cochain which value  is $1$ on
the points of $X$. The usual cap product is then defined by
$$
{\cal C}^p   (X) \otimes {\cal C}_k   (X) {\longrightarrow} {\cal
C}_{k-p} (X) \,, \quad f\otimes c   \mapsto  f\cap c
 = f\cdot c\cdot 1 = \sum _i ( -1) ^{|c_i|\cdot |f|} c_i f (c'_i)\,.
$$
The cap product with a cycle $x\in {\cal C}_k(X)$ is a well
defined homomorphism of differential graded modules, but is not a
``degree $k$ homomorphism'' of ${\cal C}^\ast (X)$-bimodules.
However,

\vspace{3mm}
 \noindent{\bf 8.3. Theorem.} {\it Let $X$ be a   path connected space and
$c \in   {\cal C}_\ast  (X) $ be a
 cycle of degree $k  >0$. Then there   exists a   (degree $k$) morphism of
${\cal C}^\ast (X)$-bimodules
$$
\gamma _c : \bar{\bB} ({\cal C}^\ast  (X),{\cal C}^\ast  (X),{\cal
C}^\ast (X)) \to   {\cal C}_\ast  (X) \,
$$
such that
\begin{enumerate}
\item[$\bullet$] $\gamma_c(1[\,]1) = c$,
\item[$\bullet$] $H(\gamma_c)\circ H(m)^{-1} : H^*(X) \to H_*(X)$ is the
cap product by $[c]$,
 $m$ is the quasi-isomorphism of ${\cal C}^*(X)$-module
defined in 1.4.
\end{enumerate}
}

\vspace{3mm} Recall that $\gamma_c$ is a degree $k$ morphism of
${\cal C}^\ast (X)$-bimodules means that the following two
properties are satisfied:

a)   $d \circ \gamma _c = (-1)^k \gamma_c \circ d$,

b) $\gamma _c (f\cdot \alpha \cdot g) = (-1) ^{|f|\,k} f\cdot
\gamma_c(\alpha)\cdot g$, \\
for $f,g \in {\cal C}^\ast  (X)$ and $ \alpha \in \bar{\bB} ({\cal
C}^\ast  (X),{\cal C}^\ast (X),{\cal C}^\ast (X))$.

\vspace{2mm}

\noindent{\bf Proof.} For   simplicity we denote by $A^e$ the
enveloping  algebra of $A={\cal C}^\ast  (X)$ and by  $ B$ the
differential graded ${\cal C}^\ast (X)$-bimodule $
\bar{\bB}({\cal C}^\ast  (X),{\cal C}^\ast (X),{\cal C}^\ast
(X))$.

Recall the loop space fibration ev : $
 X ^{S ^1}  {\to} X $, $\gamma \mapsto \gamma (0)=\gamma (1)$
with the canonical section $ \sigma : X \to X^{S^1}\,, x \mapsto$
{\small the constant loop at} $x$.  Jones defined a
quasi-isomorphism of differential  graded modules (\cite{ CJ}
Theorem 8),
$$
J_{\ast } : B\otimes _{A^e} A  \to {\cal C}^*(X^{S^1})
$$
making commutative the following    diagram of differential graded
modules
$$
\begin{array}{rccccc}
B\otimes _{A^e} A & \stackrel{ J_\ast }{\to}& {\cal C}^*(X^{S^1})\\
{\scriptstyle i}   \nwarrow & & \nearrow {\scriptstyle {\cal
C}^\ast(\mbox{\small ev}) }\\
&{\cal C}^*(X)
\end{array}
$$
where $i : {\cal C }^\ast (X) \to B \otimes _{A^e} {\cal
C}^\ast(X)\,, \quad f \mapsto 1[\,]1\otimes f $, denotes the
canonical inclusion. Let $\rho $ be the composite $ {\cal C}^\ast(
\sigma ) \circ J_*$ then $\rho$ is a retraction of $i$: $ \rho
\circ i=   id$.

Let $u\in {\cal C}^k(X)^\vee$, $k>0$, be a cycle. Using the
  canonical isomorphism of differential graded
modules
$$
\Psi : \mbox{Hom}(B \otimes _{A^e} A , \bk ) \to \mbox{Hom}_{A^e}
(B, A ^\vee)\,, \quad   \left( \Psi
 (\theta)(\alpha)\right)(f)= \theta(\alpha\otimes f)\,,
$$
we define the map
$$\theta_u :   \bar{\bB} ({\cal C}^\ast  (X),{\cal C}^\ast  (X),{\cal
C}^\ast  (X)) \to   ({\cal C}^\ast  (X) )^\vee\,, \hspace{5mm}
 \theta _u = \Psi( u\circ \rho)\,.$$ The element $\theta_u$  is
a $k$-cycle in $\mbox{Hom}_{A^e} (B, A ^\vee)$ and for any $f\in
A$, $ \theta_u(1[\,]1)(f)= u\circ \rho(1[\,]1\otimes f)= u\circ
\rho \circ i (f)= u(f)\,.$

Since the linear map
$$\lambda : {\cal C}_\ast (X) \to {\cal C}^\ast (X)^\vee$$
is   a morphism of differential graded ${\cal C}^\ast
(X)$-bimodules, for a cycle $c\in {\cal C}_k(X)$, we have a
morphism
$$\theta_{\lambda (c)} :  \bar{\bB} ({\cal C}^\ast  (X),{\cal C}^\ast
(X),{\cal C}^\ast  (X)) \to   ({\cal C}^\ast  (X) )^\vee$$ with
$\theta_{\lambda (c)}(1[\,]1) = \lambda(c)$.

Since $ \bar{\bB} ({\cal C}^\ast  (X),{\cal C}^\ast  (X),{\cal
C}^\ast (X)) $ is semifree, we deduce from  the lifting homotopy
property (cf. 1.4)
  a morphism of ${\cal C}^\ast (X)$-bimodules
$$ \gamma_c :  \bar{\bB} ({\cal C}^\ast  (X),{\cal C}^\ast  (X),{\cal
C}^\ast  (X)) \to {\cal C}_\ast (X)$$ making commutative, up to
homotopy, the   diagram
$$
\renewcommand{\arraystretch}{1.4}
\begin{array}{ccc}
 \bar{\bB} ({\cal C}^\ast  (X),{\cal C}^\ast  (X),{\cal C}^\ast  (X))
&\stackrel{\theta_{\lambda (c)}}{\longrightarrow} & {\cal
C}^\ast(X)^\vee\\
\parallel & & \uparrow {\scriptstyle \lambda }\\
 \bar{\bB} ({\cal C}^\ast  (X),{\cal C}^\ast  (X),{\cal C}^\ast  (X)) &
\stackrel{\gamma_c}{\longrightarrow} & {\cal C}_\ast (X)
\end{array}
\renewcommand{\arraystretch}{1}
$$
and such that $\gamma_c(1[\, ]1) = c$.

  The equality $H(\gamma_c)\circ H(m)^{-1} = \cap [c]$ comes from the
commutativity  of the diagram
$$
\begin{array}{ccc}
 \overline{\bB}_0({\cal C}^\ast (X) ,{\cal C}^\ast (X) ,{\cal C}^\ast (X)
) & \stackrel {\theta_{\lambda (c)} }\rightarrow &
{\cal C}^\ast (X)^\vee \\
{\scriptstyle m}\downarrow && \uparrow {\scriptstyle \lambda} \\
{\cal C}^*(X) & \stackrel{ -\cap c}{\longrightarrow} &{\cal
C}_*(X)
\end{array}
$$
i.e., for any $f,g,h \in {\cal C}^\ast(X)$, we have $ \langle
\theta_{\lambda(c)}(f[\,]g), h\rangle = \langle \lambda \circ
(-\cap c) \circ m (f[\,]g), h\rangle \,. $

\hfill $\square$

 \vspace{3mm} As a special   case, we deduce:

 \vspace{3mm}
\noindent{\bf 8.4. Theorem.} {\it  Let $M$ be a   1-connected
$\bk$-Poincar\'e duality space of  formal  dimension $d$. Then
there are quasi-isomorphisms of ${\cal C}^*(M)$-bimodules
$$
{\cal C}^\ast (M) \stackrel {m} \leftarrow \overline{\bB}({\cal
C}^\ast (M) ,{\cal C}^\ast (M) ,{\cal C}^\ast (M) )
\stackrel{\gamma}{\rightarrow} {\cal C}_\ast (M)
$$
where
 $m$ is defined in lemma 1.3
and  $\gamma= \gamma_{[M]}$ with $[M] \in H_d(M)$ a fundamental
class of $M$. In particular, the composite, $H(m) \circ
H(\gamma)^{-1} $ is the Poincar\'e isomorphism ${\cal P} : H _\ast
(M) \to H^{d-\ast}(M)$.}

 \vspace{3mm} Applying Hochschild cohomology, we obtain

 \vspace{3mm}\noindent{\bf 8.5 Theorem.} {\sl Let $M$ be a   1-connected
$\bk$-Poincar\'e duality space of  formal  dimension $d$ then
there exist   natural  linear  isomorphisms
 $$ D : H\!H^n ({\cal C}^\ast (M); {\cal C}_\ast  (M))
 \stackrel{\cong}{\longrightarrow} H\!H^{n  -d}({\cal C}^\ast  (M); {\cal
 C}^\ast  (M))\,.  $$
 }

\vspace{2mm} \noindent{\bf Proof.} Let $\varphi : N \to N'$ be a
homomorphism of differential graded $A$-bimodules and assume that
$A$ is a $\bk$-module. Then we deduce from 1.4 (see \cite{FHT-1}
for more details) that $\varphi$ induces an isomorphism of graded
modules
$$
H\!H^\ast (A;N) \to H\!H^\ast(A;N') \,.
$$
Theorem 8.5 follows directly from Theorem 8.4 when one observes
that the suspended map $s^d \gamma$ is a quasi-isomorphism of
differential graded ${\cal C}^\ast (X)$-bimodules.
\hfill{$\square$}

 \vspace{1cm}
\hspace{-1cm}\begin{minipage}{19cm} \small
\begin{tabular}{lll}
felix@math.ucl.ac.be   &jean-claude.thomas@univ-angers.fr &
vigue@math.univ-paris13.fr\\
D\'epartement de math\'ematique  &  D\'epartement de
math\'ematique  &
D\'epartement de math\'ematique \\
 Universit\'e Catholique de Louvain  &Facult\'e des Sciences  &  Institut
Galil\'ee\\
 2, Chemin du Cyclotron   &2, Boulevard Lavoisier &Universit\'e de
Paris-Nord\\
 1348 Louvain-La-Neuve, Belgium       & 49045 Angers, France & 93430
Villetaneuse, France

\end{tabular}

\end{minipage}

 \end{document}